\documentclass[11pt]{amsart}                  
\usepackage{amssymb,amsmath,amsfonts,latexsym}
\usepackage[margin=1.2in]{geometry}
\usepackage{enumerate}
\newcommand{\complex}{{\mathbb C}}

\newcommand{\reals}{{\mathbb R}}

\newcommand{\integers}{{\mathbb Z}}

\newcommand{\calg}{{\mathcal G}}
\newcommand{\calh}{{\mathcal H}}

 \DeclareMathOperator{\ind}{index}

\theoremstyle{definition}

\theoremstyle{definition}

\newcommand{\be}{\begin{eqnarray}}
\newcommand{\ee}{\end{eqnarray}}

\theoremstyle{plain}
        \newtheorem{theorem}{Theorem}[section]
        \newtheorem{lemma}[theorem]{Lemma}
        \newtheorem{remark}[theorem]{Remark}
        \newtheorem{proposition}[theorem]{Proposition}
        \newtheorem{corollary}[theorem]{Corollary}
        
\begin{document}
\title[finite dimensionality]{Hopf cyclic cohomology and Hodge theory for proper actions}
\author{Xiang Tang, Yi-Jun Yao, Weiping Zhang}
\maketitle

\begin{abstract}
We introduce a Hopf algebroid associated to a proper Lie group
action on a smooth manifold. We prove that the cyclic cohomology of
this Hopf algebroid is equal to the de Rham cohomology of invariant
differential forms. When the action is cocompact, we develop a
generalized Hodge theory for the de Rham cohomology of invariant
differential forms. We prove that every cyclic cohomology class of
the Hopf algebroid is represented by a generalized harmonic form.
This implies that the space of cyclic cohomology of the Hopf
algebroid is finite dimensional. As an application of the techniques
developed in this paper, we discuss properties of the Euler
characteristic  for a proper cocompact action.
\end{abstract}
\section{Introduction}

Let $G$ be a  Lie group, and $M$ be a smooth manifold. We assume
that $G$ acts on $M$ properly. As the $G$-action is proper, the
quotient $M/G$ is a Hausdorff stratified space. Some of the
examples of  such  spaces are  considered already in \cite{P}.

In this paper, inspired by Connes and Moscovici's Hopf cyclic theory
\cite{co-mo:gelfand-fuks}, \cite{co-mo:curved}, we introduce a Hopf
algebroid to study the ``local symmetries" of this stratified space.

Hopf algebroid was introduced by Lu \cite{lu:algebroid} in
generalizing the notion of Hopf algebra. Connes and Moscovici
\cite{connes:diffcyclic} applied this concept to generalize that of
symmetry of ``noncommutative spaces". They developed a beautiful
theory of cyclic cohomology for a Hopf algebroid, and used it to
study the transverse index theory.

Our Hopf algebroid associated to the $G$-action on $M$ is a
generalization of the Hopf algebroid introduced in the first
author's joint work with Kaminker \cite{KT}. It is shown in
\cite{KT} that if $\Gamma$ is a discrete group acting on a smooth
manifold $M$, the graded commutative algebra of differential forms
on the action groupoid $M\rtimes \Gamma$ is a topological Hopf
algebroid with the coalgebra and antipode structures defined by
taking the dual of the groupoid structure. In the case of a Lie
group $G$-action, instead of considering the algebra of differential
forms on the groupoid $M\rtimes G$, we consider the algebra
$\calh(G,M)$ of differential form valued functions on $G$. The
similar construction as in \cite{KT} defines a Hopf algebroid
structure on this algebra.

We are able to compute the cyclic cohomology of this Hopf algebroid,
which is equal to the differentiable cohomology of the groupoid
$M\rtimes G\rightrightarrows M$ with coefficient in differential
forms on $M$ considered by Crainic \cite{cr:groupoid}. As the
$G$-action is proper, Crainic's result implies that the cyclic
cohomology of the Hopf algebroid $\calh(G,M)$ is equal to the de
Rham cohomology of $G$-invariant differential forms on $M$.

Our main result of this paper is to prove a ``Hodge theorem" for
$G$-invariant differential forms on $M$ when the $G$-action is
cocompact. Our approach to this generalized Hodge theory is inspired
from the third author's joint work with Mathai \cite{MZ}. Our
strategy is to study a generalized de Rham Laplace-Beltrami operator
on the space of $G$-invariant differential forms on $M$. With some
elliptic estimates, we are able to prove that this operator has
essentially the same properties as the standard Laplace-Beltrami
operator on a compact manifold. This allows to prove that every
cyclic cohomology class of the Hopf algebroid is uniquely
represented by a harmonic form of our generalized Laplace-Beltrami
operator, which implies that the cyclic cohomology of our Hopf
algebroid is finite dimensional.

\begin{theorem}\label{thm:main} Let $G$ be a Lie group acting properly and cocompactly
on a smooth manifold $M$. The cyclic cohomology groups of
$\calh(G,M)$ are of finite dimension.
\end{theorem}

The above result allows us to introduce the Euler characteristic for
a proper cocompact action of a Lie group $G$ as the alternating sum
of the dimensions of the de Rham cohomology groups of $G$-invariant
differential forms. We are able to generalize the following two
classical results about Euler characteristic to the case of a proper
cocompact $G$-action.
\begin{enumerate}
\item  The Poincar\'e duality theorem holds for twisted de Rham
cohomology groups of $G$-invariant differential forms. In
particular, when the dimension of $M$ is odd, the Euler
characteristic of a proper cocompact $G$-action on $M$ is 0;
\item When there is a nowhere vanishing $G$-invariant vector field on $M$,
the Euler characteristic of a proper cocompact $G$-action is also 0.
\end{enumerate}

The paper is organized as follows. In Section 2, we introduce the
Hopf algebroid $\calh(G,M)$ and compute its Hopf cyclic cohomology.
In Section 3, we study a generalized Laplace-Beltrami operator, and
prove that every Hopf cyclic cohomology class of $\calh(G,M)$ can be
uniquely represented by a generalized harmonic form. In Section 4,
we introduce and study the Euler characteristic for a proper cocompact $G$-action.\\

\noindent{\bf Acknowledgments:} The  work of the first author was
partially supported by NSF grant 0703775 and 0900985. The work of
the second author was  partially supported by NSF grant 0903985. The
work of the third author was partially supported by  NNSFC and MOEC.
The first author would like to thank Marius Crainic, Niels Nowalzig,
and Hessel Posthuma for discussions about Hopf cyclic theory of Hopf
algebroids.

\section{Cyclic cohomology of Hopf algebroids}

\subsection{Hopf algebroids}

In \cite{lu:algebroid}, Lu introduced the notion of a Hopf algebroid
as a generalization of a Hopf algebra. Connes and Moscovici
\cite{connes:diffcyclic} introduced cyclic cohomology for Hopf
algebroids. Since then, many authors have studied cyclic theory for
Hopf algebroids, e.g. \cite{kh-ra:hopf-cyclic}-\cite{KP}. Oriented
by our application, we take the simplest approach for the definition
of cyclic cohomology of a Hopf algebroid, c.f.
\cite{kh-ra:hopf-cyclic}, which is also close to Connes and
Moscovici's original approach. We refer the interested readers to
\cite{ko} and \cite{KP} for the beautiful systematic study of the
general theory.

Let $A$ and $B$ be unital topological algebras. A (topological)
bialgebroid structure on $A$, over $B$, consists of the following
data.
\begin{enumerate}
\item[i)]  A continuous algebra homomorphism $\alpha:B\to A$ called the
  \emph{source map} and a continuous algebra anti-homomorphism $\beta:B\to A$ called the
\emph{target map}, satisfying $\alpha(a)\beta(b)=\beta(b)\alpha(a)$,
for all $a,b\in A$. \vskip.1in

In this paper, by tensor product $\otimes$ we always mean
topological tensor product. Let $A\otimes_B A$ be the quotient of
$A\otimes A$ by the right $A\otimes A$ ideal generated by
$\beta(a)\otimes1-1\otimes \alpha(a)$ for all   $a\in A$.

\item[ii)]   A continuous $B$-$B$
bimodule map $\Delta:
  A\to A\otimes_B A $, called the \emph{coproduct}, satisfying
    \begin{enumerate}
    \item $\Delta(1)=1\otimes 1$;
    \item $(\Delta\otimes_B Id)\Delta=(Id\otimes_B \Delta)\Delta:A\to
      A\otimes_B A\otimes_B A$,
    \item $\Delta(a)(\beta(b)\otimes1-1\otimes\alpha(b))=0$, for $a\in
      A$,  $b\in B$, 
    \item $\Delta(a_1a_2)=\Delta(a_1)\Delta(a_2)$, for $a_1, a_2\in
      A$. 
    \end{enumerate}
\item[iii)]  A continuous $B$-$B$ bimodule map $\epsilon:A\to B$, called the
  \emph{counit},  satisfying
    \begin{enumerate}
    \item $\epsilon(1)=1$;
    \item $\ker \epsilon$ is a left $ A$ ideal;
    \item $(\epsilon\otimes_B Id)\Delta=(Id\otimes_B
      \epsilon)\Delta=Id:A\to A$ 
    \item For any $a,a'\in A,\ b,b'\in B$,
    $\epsilon(\alpha(b)\beta(b')a)=b\epsilon(a)b'$, and
    $\epsilon(aa')=\epsilon(a\alpha(\epsilon(a')))=\epsilon(a\beta(\epsilon(a')))$.
    \end{enumerate}
\end{enumerate}

A topological \emph{para Hopf algebroid} is a topological
bialgebroid $A$, over $B$, which admits a continuous algebra
anti-isomorphism $S: A\to A$  such that
\[
S^2=Id,\ \ \ \ \ S\beta=\alpha,\ \ \ \ \ m_A(S\otimes_B
Id)\Delta=\beta\epsilon S:A\to A,
\]
and
\[
S(a^{(1)})^{(1)}a^{(2)}\otimes_B S(a^{(1)})^{(2)}=1\otimes_B S(a).
\]
In the above formula we have used Sweedler's notation for the
coproduct $\Delta(a)=a^{(1)}\otimes_B a^{(2)}$.

We note that in the above definition  one may allow $A$ and $B$ to
be differential graded algebras and require all of the above maps to
be compatible with the differentials and to be of degree 0. Thus one
would have a differential graded (para) Hopf algebroid (cf. \cite{gro:secondary}).

We remark that as is pointed out in \cite[Sec.2.6.13]{ko}, with our
definition any para Hopf algebroid is a Hopf algebroid as was used
in \cite{ko} and \cite{KP}. Therefore, for simplicity, in the
following, we will abbreviate ``para Hopf algebroid" to ``Hopf
algebroid".

\subsection{Hopf algebroid $\calh(G,M)$} Let $G$ be a
Lie group acting on a smooth manifold $M$.

Define $B$ to be the algebra of differential forms on $M$, and $A$
to be the algebra of $B$-valued functions on $G$. Both $A$ and $B$
are differential graded algebras with the de Rham differential. We
fix the notation that for a group element $g$ in $G$ and a smooth
function $a$ on $M$, $g^*(a)(x):=a(gx)$.

We define the source and target map $\alpha,\beta: B\to A$ as
follows,
\begin{equation}
\alpha(b)(g)=b,\qquad \text{and}\qquad \beta(b)(g)={g}^*(b).
\end{equation}
It is easy to check that $\alpha$ (and $\beta$) is an algebra (anti)
homomorphism.

When we consider the projective tensor product, the space
$A\otimes_B A$ is isomorphic to the space of $B$-valued functions on
$G\times G$, i.e.
\begin{equation}
(\phi\otimes _B \psi)(g_1, g_2)=\phi(g_1)g_1^*(\psi(g_2))
\end{equation}
for $\phi, \psi\in A$. We define the bimodule map $\Delta:
A\rightarrow A\otimes_B A$ by
\begin{equation}
\Delta(\phi)(g_1,g_2)=\phi(g_1g_2),
\end{equation}
and define the counit map $\epsilon: A\rightarrow B$ by
$\epsilon(\phi)=\phi(1)$, for $\phi\in A$.

It is straightforward to check that $(A,B, \alpha, \beta, \Delta,
\epsilon)$ is a differential graded topological bialgebroid.

To make $(A,B, \alpha, \beta, \Delta, \epsilon)$ into a Hopf
algebroid, we define the antipode on $A$ by
\begin{equation}
S(\phi)(g)={g}^*(\phi(g^{-1})).
\end{equation}

It is easy to check that $S$ satisfies properties for an antipode of
a para Hopf algebroid:

\begin{itemize}
\item $S(\beta(b))(g)=g^\ast(\beta(b)(g^{-1}))=g^\ast
(g^{-1})^\ast(b)=b$.
\item $(m_A(S\otimes Id)\Delta)(\phi)(g)= g^\ast\phi(1)$,
$(\beta\epsilon S)(\phi)(g)=g^\ast (S(\phi)(1))=g^\ast \phi(1) $.
\item One computes that $(S\otimes
Id)\Delta(a)(g_1,g_2)={g_1}^*(a((g_1)^{-1}g_2)).$ Therefore, one has
\[
\begin{split}
(\Delta S\otimes Id)\Delta(a)(g_1,g_2,g_3)&={(g_1g_2)}^*\big(a(g_2^{-1}g_1^{-1}g_3)\big),\\
S(a^{(1)})^{(1)}a^{(2)}\otimes_B
S(a^{(1)})^{(2)}(g_1,g_2)&={(g_1g_2)}^*\big(a(g_2^{-1}g_1^{-1}g_1)\big)\\
&=g_1^*\Big(g_2^*\big(a(g_2^{-1})\big)\Big)=1\otimes_B
S(a)(g_1,g_2).
\end{split}
\]
\end{itemize}

We denote this Hopf algebroid by $\calh(G,M)$.

\subsection{Cyclic Cohomology}

In this part, we briefly recall the definition of the cyclic cohomology of a
Hopf algebroid.

Let $\Lambda$ be the cyclic category. We recall the cyclic module
$A^\natural$ for $(A, B, \alpha, \beta, \Delta,\epsilon, S)$
introduced by Connes-Moscovici \cite{co-mo:curved}.

Define
\[
C^0=B,\ \ \ \ \ C^n=\underbrace{A\otimes _B A\otimes _B\cdots
\otimes_B A}_n,\ n\geq 1.
\]

Faces and degeneracy operators are defined as follows:
\[
\begin{split}
\delta_0(a^1\otimes_B\cdots \otimes_B a^{n-1})&=1\otimes_B a^1\otimes_B \cdots \otimes_B a^{n-1};\\
\delta_i(a^1\otimes_B\cdots \otimes_B a^{n-1})&=a^1\otimes_B\cdots
\otimes_B \Delta a^i\otimes_B
\cdots \otimes_B a^{n-1},\ \ \ \ \ 1\leq i\leq n-1;\\
\delta_n(a^1\otimes_B \cdots \otimes_B a^{n-1})&=a^1\otimes_B\cdots \otimes_Ba^{n-1}\otimes_B1;\\
\sigma_i(a^1\otimes_B\cdots \otimes_B a^{n+1})&=a^1\otimes_B\cdots
\otimes_B a^i\otimes_B \epsilon(a^{i+1})\otimes_Ba^{i+2}\otimes_B
\cdots \otimes_B a^{n+1}.
\end{split}
\]

The cyclic operators are given by
\[
\tau_n(a^1\otimes_B\cdots \otimes_B
a^n)=(\Delta^{n-1}S(a^1))(a^2\otimes \cdots a^n\otimes 1).
\]
The cyclic cohomology of $(A, B, \alpha, \beta, \Delta, \epsilon,
S)$ is defined to be the cyclic cohomology of $A^\natural$.

\subsection{Hopf cyclic cohomology of the Hopf algebroid $\calh(G,M)$}
In this section we explain the  computation the Hopf cyclic
cohomology of the Hopf algebroid $\calh(G,M)$.

We review briefly the definition of differentiable cohomology of a
Lie group. Let $G$ be a Lie group acting on a manifold $M$. Consider
$E$ a $G$-equivariant bundle on $M$. An $E$-valued differentiable
$p$-cochain is a smooth map $c$ mapping $G^{\times p}$ to a smooth
section of $E$, i.e. $C_d^p(G;E)=C^\infty(G^{\times p}; \Gamma(E))$.
The differential $d$ on $C^\bullet_d(G;E)$ is defined by
\begin{eqnarray*}
(dc)(g_1, \cdots, g_{p+1})&=&g_1^*(c(g_2, \cdots, g_{p+1}))+
\sum_{i=1}^p(-1)^ic(g_1, \cdots, g_ig_{i+1}, \cdots,
g_{p+1})\\
&&\hspace{1.5cm}+(-1)^{p+1}c(g_1, \cdots, g_p).
\end{eqnarray*}
The differentiable cohomology $H_d^\bullet(G;E)$ of $G$ with
coefficient $E$ is defined to be the cohomology of
$(C^\bullet_d(G;E),d)$. We remark that the space $C_d^\bullet(G;E)$
has a structure of a cyclic simplicial space. We recall its
definition below,
\[
\delta_i(a)( g_1, \cdots, g_n,
g_{n+1})=\left\{\begin{array}{ll}g_1^*(a(g_2, \cdots, g_{n+1}))& i=0\\
a(  g_1, \cdots, g_ig_{i+1}, \cdots, g_{n+1})& 1\leq i\leq n\\
a(g_1, \cdots, g_{n})& i=n+1\end{array}\right.,
\]
and
\[
\sigma_i(a)( g_1, \cdots, g_n)=a(  \cdots, g_{i-1}, 1, g_{i},
\cdots, g_n),
\]
and
\[
t(a)( g_1, \cdots, g_n)=(g_1\cdots g_n)^*a((g_1g_2\cdots g_n)^{-1},
g_1, \cdots, g_{n-1}).
\]
The cohomology of this simplicial complex $C_d^\bullet(G,E)$ is isomorphic
to the differentiable cohomology of $G$ with coefficient in $E$. We
are now ready to present the computation of Hopf cyclic cohomology
of the Hopf algebroid $\calh(G,M)$.

\begin{theorem}\label{thm:cohomology-hopf}
Let $G$ be a Lie group acting on a smooth manifold. We have
\[
HC^\bullet(\calh(G,M))=\oplus_{k\geq
0}H^{\bullet-2k}(G;(\Omega^*(M),d)).
\]
\end{theorem}
\noindent{\bf Proof:} We observe that a $p$-cochain on $\calh(G, M)$
can be identified with $\Omega^*(M)$-valued functions on $G^{\times
p}$, $p\geq 0$. This identification respects the cyclic simplicial
structures on $C^\infty(G^{\times \bullet}, \Omega^*(M)) $
and $\calh(\calg)^{\natural}$. Therefore, we conclude that the
Hochschild cohomology of $\calh(G,M)$ is isomorphic to the
differentiable cohomology $H^\bullet(G;(\Omega^*(M),d))$. By the
SBI-sequence of cyclic cohomology, we have
\[
HC^\bullet(\calh(G,M))=\oplus_{k\geq
0}H^{\bullet-2k}(G;(\Omega^*(M),d)).\qquad \qquad \Box
\]

\medskip

Let $\Omega^*(M)^G$ be the space of $G$-invariant differential forms
on $M$, which inherits a natural de Rham differential $d$. By
\cite[Section 2.1, Prop 1]{cr:groupoid}, if $G$ acts on $M$
properly, then the differentiable cohomology
$H^\bullet(G;(\Omega^*(M),d))$ is computed as follows,
\[
H^\bullet(G;(\Omega^*(M),d))=H^\bullet(\Omega^*(M)^G,d).
\]
\begin{proposition}\label{prop:inv-hom}If $G$ acts on $M$ properly, then we have
\[
HC^\bullet(\calh(G,M))=\oplus_{k\geq
0}H^{\bullet-2k}(\Omega^*(M)^G,d)\] and
\[ HP^\bullet(\calh(G,M))=\oplus_{k\in
\integers}H^{\bullet+2k}(\Omega^*(M)^G,d).
\]
\end{proposition}

\section{Generalized Hodge theory and a Proof of Theorem 1.1}
\label{sec:hodge} Now we prove Theorem 1.1. According to Proposition
\ref{prop:inv-hom}, all we need to prove is:
\smallskip

\centerline{ The cohomology groups $H^{\bullet}(\Omega^*(M)^G,d)$
are of finite dimension. }
\medskip

We adapt the proof of the finite dimensionality of the de Rham
cohomology of compact manifolds. (cf. \cite[Ch. 6]{W})

Without loss of generality, as $G$ acts on $M$ properly, we may
assume that $M$ is endowed with a $G$-invariant metric. And since
$M/G$ is compact, there exists a compact subset $Y$ of $M$ such that
$\displaystyle G(Y)=\bigcup_{g\in G} gY=M$ (cf. \cite[Lemma
2.3]{P}). With the usual de Rham Hodge $*$ operator, on $\Omega^*(M)$ we
consider the following inner product
\begin{equation}\label{2.5}
\langle \alpha, \beta\rangle_0=\int_M \alpha\wedge *\beta.
\end{equation}

As $Y$ is a closed subset of $M$, there exist $U$, $U'$, two open
subsets of $M$, such that $Y\subset U$ and that the closures
$\overline{U}$ and $\overline{U'}$ are both compact in $M$, and that
$\overline{U}\subset U'$. It is easy to construct a smooth function
$f:M\rightarrow [0,1]$ such that $f|_U=1$ and ${\rm Supp}(f)\subset
U'$.

Let $\Gamma(\Omega^*(M))^G$ be the subspace of $G$-invariant
sections of $\Omega^*(M)$. For an open set $W$ of $M$, define
\begin{align}\label{2.7}
\|s\|_{W,\,0}^2=\int_W\left\langle s(x),s'(x)\right\rangle dx,\qquad
\|s\|_{W,\,1}^2=\|s\|_{W,\,0}+\langle\Delta(s),s\rangle_{W,0}.
\end{align}

For any $s\in\Gamma(\Omega^*(M))^G$,
\begin{align}\label{2.6.1}
\|s\|_{U,\,0}\leq \|fs\|_0\leq \|s\|_{U',\,0}.
\end{align}
As $G(Y)=M$, $G(U)=M$. Since $\overline{U'}$ is compact,  there are finitely many elements $g_1,
\cdots, g_k$ of $G$ such that $g_1U\cup \cdots \cup g_kU$ covers
$U'$.  If $s$ is a $G$-invariant section of $\Omega^*(M)$, it is
easy to see that there exists a positive constant $C>0$,
\begin{align}\label{2.6}
\|s\|_{U',\,0}\leq C\|s\|_{U,\,0}.
\end{align}

Let $dg:=dm(g)$ be the right invariant Haar measure on $G$. Define
$\chi:G\to \reals^+$ by $dm(g^{-1})=\chi(g) dm(g)$. We define ${\bf
H}^0_f(M,\Omega^*(M))^G$ to be the completion of the space $\{fs:
s\in \Gamma(\Omega^*(M))^G\}$ under the norm $\|\cdot\|_0$
associated to the inner product (\ref{2.5}). As indicated in in the
Appendix in \cite{MZ} written by Bunke, we will first prove the
following proposition.

\begin{proposition}\label{prop:projection}For any $\mu\in L^2(M,\Omega^\ast(M))$, define
\begin{equation}\label{Pf}
\left(P_f \mu\right)(x) = \frac{f(x)}{(A(x))^2}\int_G \chi(g) f(gx)
\mu(gx) \, dg,
\end{equation}
where
\begin{equation}\label{Ax}
A(x)=\left(\int_G \chi(g)(f(gx))^2 \, dg \right)^{1/2}
\end{equation}
is a $G$-equivariant function on $M$, i.e.
$A(gx)^2=\chi(g)^{-1}A(x)^2$, and is strictly positive. The operator
$P_f$ defines an orthogonal projection from $L^2(M,\Omega^\ast(M))$
onto ${\bf H}^0_f(M,\Omega^*(M))^G$.
\end{proposition}

\begin{proof}
It is straightforward to check that $A(x)$ is strictly positive and
equivariant. In order to show that (\ref{Pf}) defines actually an
orthogonal projection, we need to prove the following properties.
\begin{itemize}
\item $P_f^2=P_f$. For any $\mu\in L^2(M,\Omega^\ast(M))$,
\begin{eqnarray}
(P_f^2 \mu)(x) & = & \frac{f(x)}{(A(x))^2}\int_G \chi(g) f(gx)
\left(P_f \mu\right)(gx) \, dg\cr &=& \frac{f(x)}{(A(x))^2}\int_G
\chi(g) f(gx) \left( \frac{f(gx)}{(A(gx))^2}\int_G \chi(h) f(hgx)
\mu(hgx) \, d h \right) \, dg\cr &=& \frac{f(x)}{(A(x))^2}\int_G
\frac{\chi(g)(f(gx))^2}{(A(x))^2} \, dg \left( \int_G \chi(hg)
f(hgx) \mu(hgx) \, d (hg) \right) \cr &=& \frac{f(x)}{(A(x))^2}
\left( \int_G \chi(hg) f(hgx) \mu(hgx) \, d (hg) \right) \cr & = &
(P_f \mu)(x).
\end{eqnarray}
\item $P_f$ is self-adjoint. For any $\mu, \nu \in L^2(M,\Omega^\ast(M))$,
\begin{eqnarray}
\langle P_f \mu, \,\nu \rangle_0 &=& \int_M
\frac{f(x)}{(A(x))^2}\int_G \chi(g) f(gx) \langle \mu(gx),\,
\nu(x)\rangle \, dg\, dx\cr &=& \int_M \int_G\frac{f(g^{-1}
x^\prime)}{(A(g^{-1}x^\prime))^2} \chi(g) f(x^\prime) \langle
\mu(x^\prime),\, \nu(g^{-1}x^\prime)\rangle \, dg\, dx^\prime
\hspace{1.5cm} (x^\prime=gx)\cr &=& \int_M\int_G
\frac{f(x^\prime)}{(A(x^\prime))^2} f(g^{-1} x^\prime) \langle
\mu(x^\prime),\, \nu(g^{-1}x^\prime)\rangle \, dg\, dx^\prime \cr
&=& \int_M\int_G \frac{f(x^\prime)}{(A(x^\prime))^2} f(g^{-1}
x^\prime) \langle \mu(x^\prime),\, \nu(g^{-1}x^\prime)\rangle \,
\chi(g^{-1})d(g^{-1})\, dx^\prime\cr &=& \langle
 \mu, \,P_f \nu \rangle_0.
\end{eqnarray}
\item it is also straightforward to check that for $\mu=f\alpha \in {\bf
H}^0_f(M,\Omega^*(M))^G$ where $\alpha\in \Omega^*(M)^G$, $P_f
\mu=\mu$.
\end{itemize}
The Proposition is thus proved.
\end{proof}

Define ${\bf H}^1_f(M,\Omega^*(M))^G$ to be the completion of $\{fs:
s\in \Gamma(\Omega^*(M))^G\}$ under a (fixed) first Sobolev norm
associated to the inner product (\ref{2.5}). And in general define
${\bf H}^{k}_f(M,\Omega^*(M))^G$ (and ${\bf
H}^{-1}_f(M,\Omega^*(M))^G$) to be the completion of the space
$\{fs: s\in \Gamma(\Omega^*(M))^G\}$ under the corresponding ${\bf
H}_f^{k}$ (and ${\bf H}_f^{-1}$ norm) for $k\geq 2$. (For any open
subset $W$ of $M$, and any compactly supported smooth differential
form $s$ on $W$, $\|s\|_{W,\, k}^2=\|(1+\Delta)^{k/2}(s)\|_{W,\,0}^2$,
for $k\geq 2$.)

This time we investigate the operator
\begin{eqnarray}\label{eq:dfn-df}
d_f : {\bf H}_f^{1} & \rightarrow & {\bf H}_f^{0}\cr
 f\alpha &\mapsto & f \, d\alpha.
\end{eqnarray}
We also consider its adjoint $d_f^\ast: {\bf H}_f^{1}  \rightarrow
{\bf H}_f^{0}$: for $\alpha\in \Gamma(\Omega^p(M))^G, \beta \in \Gamma(\Omega^{p+1}(M))^G$,
\begin{equation}
\langle d_f(f\alpha),\, f\beta\rangle_0=\langle f\alpha,\, d_f^*(f\beta)\rangle_0.
\end{equation}

For $f\in C^\infty(M)$, denote $\nabla(f)$ to be the gradient vector
field associated to the riemannian metric on $M$. We can show easily
that
\begin{equation}\label{dfstar}
d_f^\ast (f\beta)= P_f(-2 i_{\nabla f} \beta)+ f\delta\beta,
\end{equation}
where $\delta= (-1)^{n(p+1)+n+1}* d *: \Omega^{p+1}(M)\rightarrow \Omega^p(M)$ and $i_V \alpha$ is the contraction of
the form $\alpha$ with the vector field $V$.

Now we define a self-adjoint operator
\begin{equation}\label{Deltatilde}
\tilde{\Delta}= d_f d_f^\ast+ d_f^\ast d_f.
\end{equation}

\begin{proposition}\label{prop:fredholm}
$\tilde{\Delta}: {\bf H}_f^{2}  \rightarrow  {\bf H}_f^{0}$ is
Fredholm.
\end{proposition}

\begin{proof}
We will prove this fact by establishing a G\aa rding type inequality. Let $f\alpha\in {\bf H}^{0} $,
\begin{eqnarray}\label{19}
\tilde{\Delta}(f\alpha) & = & d_f (d_f^\ast(f\alpha))+ d_f^\ast(
d_f(f\alpha))\cr &=& d_f (-P_f(2 i_{\nabla f} \alpha)+
f\delta\alpha) + d_f^\ast(f d \alpha)\cr &=&- d_f (P_f(2 i_{\nabla
f} \alpha)) + f d\delta \alpha- P_f(2 i_{\nabla f} d \alpha)+f\delta
d\alpha\cr &=& f\Delta \alpha - d_f (P_f(2 i_{\nabla f} \alpha)) -
P_f(2 i_{\nabla f} d \alpha).
\end{eqnarray}

Now using (\ref{Pf}), we have the following estimates:

\begin{eqnarray}
&&\left\| d_f \left( P_f(2 i_{\nabla f} \alpha) \right)\right\|_0\cr
&=& \left\| 2 f(x) d\left(\frac{\int_G \chi(g) f(gx) (i_{\nabla f}
\alpha)(gx) \, dg }{(A(x))^2}\right) \right\|_0\\
&=& 2\left\|f
\int_G \chi(g) d\left( \frac{g^\ast f}{A^2}\right) ( g^\ast
i_{\nabla f} \alpha) \, dg  \ + f\int_G \chi(g) \frac{g^\ast f}{A^2}
d\left( g^\ast i_{\nabla f} \alpha \right)\, dg \right\|_{0}\cr
&\leqslant & 2 \left\| f\int_G \chi(g) d\left( \frac{g^\ast
f}{A^2}\right) ( g^\ast i_{\nabla f} \alpha) \, dg\right\|_{0}\  \cr
& & \hspace{3cm} + 2\left\| f\int_G \chi(g) \frac{g^\ast f}{A^2}
d\left( g^\ast i_{\nabla f} \alpha \right)\, dg \right\|_{0},\cr
\end{eqnarray}
where

\begin{eqnarray}
\left\|f \int_G \chi(g) d\left( \frac{g^\ast f}{A^2}\right) ( g^\ast
i_{\nabla f} \alpha) \, dg\right\|_{0} &\leqslant & \left\|f \int_G
\chi(g)  \left| d\left( \frac{g^\ast f}{A^2}\right)\right| \left|(
g^\ast i_{\nabla f} \alpha)\right| \, dg\right\|_{0}\cr &\leqslant &
\left\|f \int_G \chi(g) \left| d\left( \frac{g^\ast
f}{A^2}\right)\right| \left| g^\ast (\nabla f)\right| \left|
g^\ast\alpha\right| \, dg\right\|_{0}\cr &\leqslant & \left\|f
\int_G \chi(g) \left| d\left( \frac{g^\ast f}{A^2}\right)\right|
\left| g^\ast (\nabla f)\right| \, dg  \left|
\alpha(x)\right|\right\|_{0}.
\end{eqnarray}
As $f$ has a compact support,
\[f\int_G
\chi(g)  \left| d\left( \frac{g^\ast f}{A^2}\right)\right| \left|
g^\ast (\nabla f)\right|  \, dg \] is finite and continuous
everywhere,  and therefore is bounded from above by a constant on
$\overline{U}$. Hence we have
\begin{equation}
 2\left\|f \int_G
\chi(g) d\left( \frac{g^\ast f}{A^2}\right) ( g^\ast i_{\nabla f}
\alpha) \, dg\right\|_{0} \leqslant  \tilde{C}_1
\|\alpha\|_{U',\,0}\leqslant \tilde{C}_2\|f\alpha\|_{0}.
\end{equation}

Consider
\begin{equation}\label{pfdialpha}
2\left\| f\int_G \chi(g)  \frac{g^\ast f}{A^2} d\left( g^\ast
i_{\nabla f} \alpha \right)\, dg  \right\|_{ 0} =2\left\|P_f(d\circ
i_{\nabla(f)}\alpha)\right\|_{0}\leq 2\|d\circ
i_{\nabla(f)}\alpha\|_{0}.
\end{equation}

Notice that $f$ is supported inside $U'$. We choose a cut-off
function $c$ which is 1 on the support of $f$ and 0 outside $U''
\subset \overline{U''}\subset U'$. We have
\begin{equation}\label{dicalpha}
\|d\circ i_{\nabla(f)}\alpha\|_{0}=\|d\circ i_{\nabla(
f)}(c\alpha)\|_{U',\,0}.
\end{equation}

We observe that $d\circ i_{\nabla(f)}$ is a differential operator on
$\Omega^*(M)$ of order 1. As $c\alpha$ is a compactly supported
smooth function in $U'$, we have
\begin{equation}\label{dialpha}
\|d\circ i_{\nabla( f)}(c\alpha)\|_{U',\,0}\leq C
\|c\alpha\|_{U',\,1}.
\end{equation}

We compute $\|c\alpha\|_{U',\,1}$ to be
\begin{eqnarray}\label{calpha}
\|c\alpha\|^2_{U',\,1}&=&\|c\alpha\|^2_{U',\,0}+<\Delta(c\alpha),c\alpha>_{U'}\cr
&=&\|c\alpha\|^2_{U',\,0}+<d(c\alpha),d(c\alpha)>_{U'}+<\delta(c\alpha),
\delta(c\alpha)>_{U'}\cr
&=&\|c\alpha\|^2_{U',\,0}+\|d(c\alpha)\|_{U',\,0}+\|\delta(c\alpha)\|_{U',\,0}.
\end{eqnarray}
We discuss one by one the terms in the above equation.

As $c$ is bounded from above by 1, we have
\[
\|c\alpha\|_{U',\,0}\leq \|\alpha\|_{U',\,0}.
\]

As $c$ and $dc$ are both bounded,
\begin{equation}\label{dcalpha}
\|d(c\alpha)\|_{U',\,0}=\|dc\wedge \alpha+cd\alpha\|_{U',\,0}\leq
\|dc\wedge \alpha\|_{U',\, 0}+\|cd\alpha\|_{U',\,0}\leq
C_1\|\alpha\|_{U',\,0}+C_2\|d\alpha\|_{U',\,0}.
\end{equation}

Similarly, as $c$ and $\nabla(c)$ are compactly supported, they are
both bounded. We have
\begin{equation}\label{deltacalpha}
\|\delta(c\alpha)\|_{U',\,0}=\|i_{\nabla(c)}\alpha+c\delta\alpha\|_{U',\,0}\leq
C_3\|\alpha\|_{U',\,0}+\|\delta \alpha\|_{U',\,0}.
\end{equation}

As $G(U)=M$ and $\overline{U'}$ is compact, there are finitely many
$g_1, \cdots, g_k$ such that $U'\subset g_1U\cup \cdots \cup g_kU$.
As $\alpha$, $d\alpha$, and $\delta\alpha$ are all $G$-invariant, we
have
\begin{eqnarray}\label{uu'}
\|\alpha\|_{U',\,0}  &\leq  C_4\|\alpha\|_U& \leq
C_4\|f\alpha\|_{0}\cr \|d\alpha\|_{U',\,0} &\leq  C_5\|d\alpha\|_U&
\leq C_5\|d(f\alpha)\|_{0}\cr \|\delta \alpha\|_{U',\,0}&\leq
C_6\|\delta \alpha\|_U&\leq C_6\|\delta (f\alpha)\|_{0}.
\end{eqnarray}

Summarizing inequalities (\ref{pfdialpha})-(\ref{uu'}), we have that
\[
2\left\|f \int_G \chi(g)  \frac{g^\ast f}{A^2} d\left( g^\ast
i_{\nabla f} \alpha \right)\, dg  \right\|_{ 0}\leq A\|f\alpha\|_1.
\]

Similarly,

\begin{equation}
\left\| P_f(2 i_{\nabla f} d\alpha)  \right\|_0 \leqslant  2\left\|
i_{\nabla f} d\alpha \right\|_0 \leqslant \tilde{C}_7  \left\|
d\alpha \right\|_{U^\prime,\, 0}  \leqslant \tilde{C}_8  \left\|
d(f\alpha) \right\|_{0}\leq B\|f\alpha\|_1.
\end{equation}

By combining these inequalities, we have

\begin{eqnarray}\label{tildedeltafalpha}
\|\tilde{\Delta}(f\alpha)\|_0 &\geqslant & \| f\Delta \alpha\|_0 -
\| d_f (P_f(2 i_{\nabla f} \alpha))\|_0 - \| P_f(2 i_{\nabla f} d
\alpha)\|_0 \cr &\geqslant & \| f\Delta \alpha\|_0 -
(A+B)\|f\alpha\|_1\cr &\geqslant &\left\|\Delta\alpha \right\|_{U,\,
0} - \tilde{B}\|f\alpha\|_1 .
\end{eqnarray}

By the standard elliptic inequality, we have that
\begin{equation}\label{deltaalpha0}
\|\Delta\alpha\|_{U,0}\geq
\tilde{A}\|\alpha\|_{U,\,2}-D\|\alpha\|_{U,\,0}.
\end{equation}

By definition, we have
\begin{equation}\label{falpha2}
\|f\alpha\|_{2}^2=\|f\alpha\|_0^2+2\|(d+\delta)(f\alpha)\|_0^2+\|(d+\delta)^2(f\alpha)\|^2_0.
\end{equation}

We compute
\begin{eqnarray}
(d+\delta)^2(f\alpha)&=&(d+\delta)(df\wedge \alpha +fd\alpha
-i_{\nabla(f)}\alpha+f\delta(\alpha))\cr &=&\delta(df\wedge
\alpha)-i_{\nabla (f)}d\alpha+f\delta
d\alpha-d(i_{\nabla(f)}\alpha)+df\wedge\delta(\alpha)+fd\delta
\alpha\cr &=& f\Delta(\alpha)+\delta(df\wedge \alpha)-i_{\nabla
(f)}d\alpha-d(i_{\nabla(f)}\alpha)+df\wedge\delta(\alpha).
\end{eqnarray}

We notice that $\delta(df\wedge \alpha)$, $i_{\nabla (f)}d\alpha$,
$d(i_{\nabla(f)}\alpha)$, and $df\wedge\delta(\alpha)$ are all
differential operators on $\alpha$ of order less than or equal to 1.
So similar estimates as (\ref{dicalpha})-(\ref{uu'}) show that every
piece of them is bounded by a multiple of $\|f\alpha\|_1$.

We prove by the similar arguments as (\ref{uu'}) that
\begin{equation}\label{deltaalphauu'}
\|f\Delta\alpha\|_0\leq \|\Delta\alpha\|_{U',\,0}\leq D_1\|\Delta
\alpha\|_{U,\,0}\leq D_1\|\alpha\|_{U,\,2}.
\end{equation}

By (\ref{falpha2})-(\ref{deltaalphauu'}), we have
\[
\|\alpha\|_{U,\,2}\geq D_2\|f\alpha\|_2-D_3\|f\alpha\|_1.
\]

With the above estimate, from (\ref{tildedeltafalpha}) and
(\ref{deltaalpha0}), we have
\[
\|\tilde{\Delta}(f\alpha)\|_0\geq D_3\|f\alpha\|_2-D_4\|f\alpha\|_1.
\]

By the so-called Peter-Paul Inequality, we have
\begin{equation}
 \|f\alpha\|_1\leqslant \frac{1}{2}D_3/D_4 \|f\alpha\|_2+ D_5 \|f\alpha\|_0.
\end{equation}

In summary, we have
\begin{equation}\label{garding}
\|\tilde{\Delta}(f\alpha)\|_0 \geqslant \frac{1}{2}D_3 \|f\alpha\|_2
-D_6 \|f\alpha\|_0.
\end{equation}
Due to the fact that the embedding of ${\bf
H}^{2}_f(M,\Omega^*(M))^G$ in ${\bf
H}^{0}_f(M,\Omega^*(M))^G$ is compact, the above G\aa rding type inequality implies that
 $\tilde{\Delta}$ is Fredholm.
\end{proof}

\begin{corollary}\label{cor:dimker}
$\dim (\mathrm{ker}\, \tilde{\Delta})=\dim(\mathrm{coker}\,
\tilde{\Delta})<+\infty $.
\end{corollary}

\begin{lemma}
$\mathrm{ker}\, \tilde{\Delta}=(\mathrm{Im}\, \tilde{\Delta})^\bot \cap {\bf
H}^{2}_f$.
\end{lemma}

\begin{proof}
We have
\begin{eqnarray*}
 & & f\alpha \in (\mathrm{Im}\, \tilde{\Delta})^\bot \cap {\bf
H}^{2}_f \\
&\Leftrightarrow & \forall f\beta\in {\bf
H}^{2}_f,\,\,\langle \tilde{\Delta}(f\beta),\,f\alpha \rangle_0=0,\\
&\Leftrightarrow & \forall f\beta\in {\bf
H}^{2}_f,\,\,\langle f\beta,\, \tilde{\Delta}(f\alpha)\rangle_0=0.
\end{eqnarray*}
As ${\bf
H}^{2}_f$ is dense in ${\bf
H}^{0}_f$, so $\tilde{\Delta}(f\alpha)=0$, which  is equivalent to $f\alpha\in \mathrm{ker}\, \tilde{\Delta}$.

\end{proof}

This lemma together with the previous Corollary \ref{cor:dimker}
implies that
\begin{equation}
\mathrm{ker}\, \tilde{\Delta}=(\mathrm{Im}\, \tilde{\Delta})^\bot,
\end{equation}
i.e., we have the decomposition
\begin{equation}\label{HodgeDecomp}
{\bf
H}^{0}_f = \mathrm{ker}\, \tilde{\Delta}\oplus \mathrm{Im}\, \tilde{\Delta}.
\end{equation}
Therefore, we can define the projection $H : {\bf
H}^{0}_f  \rightarrow \mathrm{ker}\, \tilde{\Delta}$. Let $f\alpha \in {\bf
H}^{0}_f$, then $f\alpha - H(f\alpha)\in  \mathrm{Im}\, \tilde{\Delta}$. So there is a unique $f\beta \in  \mathrm{Im}\, \tilde{\Delta}$ such that
\[\tilde{\Delta}(f\beta)=f\alpha - H(f\alpha).\]
We define in this way the Green operator $\mathfrak {G}:
f\alpha\mapsto f\beta$.

We will need the following propositions to explore the properties of the Green operator.

\begin{proposition}\label{6.6}
Let $\{f\alpha_n\}$ be a sequence of smooth $p$-forms in ${\bf
H}^{2}_f(M,\Omega^*(M))^G$ such that $\|f\alpha_n\|_0\leqslant c$
and $\|\tilde{\Delta} (f\alpha_n)\|_0\leqslant c$ for all $n$ and
for some constant $c>0$. Then it has a Cauchy subsequence.
\end{proposition}

\begin{proof}

We prove that $f\alpha_n$ is a bounded sequence in ${\bf H}^1_f(M,
\Omega^*(M))^G$. Then we conclude the proposition by the fact that
${\bf H}^1_f(M, \Omega^*(M))^G$ is compactly embedded in ${\bf
H}^0_f(M, \Omega^*(M))^G$.

We have the following estimates.
\[
\|f\alpha_n\|^2_{1}=\|f\alpha_n\|^2_0+\langle \Delta(f\alpha_n),
f\alpha_n\rangle_{0}.
\]

By the Cauchy-Schwarz inequality, we have
\[
\langle \Delta(f\alpha_n), f\alpha_n\rangle_{0}\leq
\|\Delta(f\alpha_n)\|_0\|f\alpha_n\|_0.
\]

By Inequality (\ref{garding}), we have that
\[
\|\Delta(f\alpha_n)\|_0\leq \|f\alpha_n\|_2\leq
A\|\tilde{\Delta}(f\alpha_n)\|_0+B\|f\alpha_n\|_0\leq (A+B)c.
\]

Therefore, $\|f\alpha_n\|_1$ is bounded by $c\sqrt{A+B+1}$.
\end{proof}

Now we prove the regularity for $\tilde{\Delta}$:

\begin{proposition}\label{prop:reg}
 If $f\beta$ is ${\bf H}_f^k(M, \Omega^*(M))^G$ and
 \begin{equation}
\tilde{\Delta}(f\alpha) = f\beta
 \end{equation}
on $M$, then $f\alpha$ is ${\bf H}_f^{k+2}(M, \Omega^*(M))^
G$ for any $k\geq 0$. In particular, if $f\beta$ is a smooth differential form, so is $f\alpha$.

\end{proposition}

\begin{proof} As $f$ is smooth and compactly supported, it is
sufficient to prove the differentiability of $\alpha$. This is a
local statement. As $\alpha$ and $\beta$ are both $G$-invariant and
$G(U)=M$, we can restrict our analysis to $U$.

By (\ref{19}) and (\ref{Pf}),

\begin{eqnarray}
\tilde{\Delta}(f\alpha) & = & f\Delta \alpha - d_f (P_f(2 i_{\nabla
f} \alpha)) - P_f(2 i_{\nabla f} d \alpha)\cr &=& f\Delta \alpha -
d_f \left( 2 f \int_G \chi(g)\frac{g^\ast f}{A^2} g^\ast
\left(i_{\nabla f} \alpha \right)\, dg \right)\cr &&\hspace{1.5cm}-
2 f \int_G \chi(g) \frac{g^\ast f}{A^2} g^\ast \left(i_{\nabla f}
d\alpha \right)\, dg.
\end{eqnarray}
As $\alpha$ is $G$-invariant (so is
$d\alpha$), we can find two $G$-invariant smooth vector
fields $V_1$, $V_2$ (which depend only on $f$) such that
\begin{equation}\label{reg}
\tilde{\Delta}(f\alpha)=  f\Delta \alpha + f d(i_{V_1}\alpha)
+ f i_{V_2} d\alpha.
\end{equation}

Notice that on $U$, $f=1$. Hence Equation (\ref{reg}) implies that
on $U$
\begin{equation}
\beta = \Delta \alpha +  d(i_{V_1}\alpha) +  i_{V_2} d\alpha.
\end{equation}
The last two terms are of lower order, so the regularity of
$\tilde{\Delta}$ is implied by that of $\Delta$, the usual
Laplace-Beltrami operator on $M$.
\end{proof}

The Green operator $\mathfrak {G}$ has the following properties:
\begin{enumerate}
\item[$1^\circ$] $\mathfrak {G}$ is bounded. To this end we need to prove the existence of a constant $c>0$ such that  for any $f\beta\in \mathrm{Im}\, \tilde{\Delta}$,
\begin{equation}\label{fbeta}
\|f\beta\|_0\leqslant c \| \tilde{\Delta}(f\beta) \|_0.
\end{equation}

Suppose the contrary, then there exists a sequence $f\beta_j \in
\mathrm{Im}\, \tilde{\Delta} $ with
\begin{equation}
\| f\beta_j \|_0=1\,\,\,\text{and}\,\,\, \|\tilde{\Delta}(f\beta_j
)\|_0\rightarrow 0.
\end{equation}
By Proposition \ref{6.6}, $\{f\beta_j\}$ has a Cauchy subsequence,
which one can assume to be $\{f\beta_j\}$ itself without loss of
generality. Hence $\displaystyle \lim_{j\rightarrow \infty} \langle
f\beta_j,\, f\psi\rangle_0 $ exists for each $f\psi\in {\bf
H}^{0}_f(M,\Omega^*(M))^G$. It defines a linear functional $l$ which
is clearly bounded, and

\begin{equation}\label{weaksol}
l(\tilde{\Delta}(f\varphi ))=\lim_{j\rightarrow \infty} \langle
f\beta_j,\, \tilde{\Delta}(f\varphi )\rangle_0 = \lim_{j\rightarrow
\infty} \langle \tilde{\Delta}(f\beta_j),\, f\varphi \rangle_0=0.
\end{equation}
We obtain the existence of $f\beta \in (\mathrm{Im}\,
\tilde{\Delta})^\bot= \ker{\tilde{\Delta}}$ such that
\begin{equation}
l(f\psi)=\langle f\beta,\, f\psi\rangle_0,\qquad \text{with}\quad
f\beta_j\rightarrow f\beta \quad \text{in}\quad {\bf H}_f^0(M,
\Omega^*(M))^G.
\end{equation}

\noindent{From} Equation (\ref{weaksol}), we know that $f\beta$ is a weak
solution of $\tilde{\Delta}(\xi)=0$. It follows from Proposition
\ref{prop:reg} that $f\beta$ is actually smooth and a strong
solution of $\tilde{\Delta}(\xi)=0$. Now as $\|f\beta_j\|_0=1$ and
$f\beta_j \in \mathrm{Im}\, \tilde{\Delta} $, it follows that
$\|f\beta\|_0=1$ and $f\beta \in \mathrm{Im}\, \tilde{\Delta} $.
Hence, $f\beta\in \mathrm{Im}\, \tilde{\Delta}\cap \mathrm{Im}\,
\tilde{\Delta}^\perp=\{0\}$, which yields a contradiction.

\item[$2^\circ$] $\mathfrak {G}$ is self-adjoint. In fact,
\begin{eqnarray}
\langle \mathfrak {G}(f\alpha),\, f\beta\rangle_0 &=& \langle
\mathfrak{G}(f\alpha),\, f\beta - H(f\beta)\rangle_0 = \langle
\mathfrak{G}(f\alpha),\, \tilde{\Delta}(\mathfrak
{G}(f\beta))\rangle_0 \cr &=&\langle \tilde{\Delta}(\mathfrak
{G}(f\alpha)),\, \mathfrak {G}(f\beta)\rangle_0 = \langle
f\alpha-H(f\alpha),\, \mathfrak {G}(f\beta)\rangle_0 =\langle
f\alpha,\, \mathfrak {G}(f\beta)\rangle_0.
\end{eqnarray}
\item [$3^\circ$] $G$ maps a bounded sequence into one with
Cauchy subsequences, due to the fact that the embedding of ${\bf
H}^{2}_f$ into ${\bf H}^{0}_f$ is compact.
\end{enumerate}

Moreover, we have

\begin{proposition}
The Green operator $\mathfrak {G}$ commutes with any linear operator
that commutes with $\tilde{\Delta}$.
\end{proposition}

\begin{proof}
Suppose that $T: f\Gamma(M,\Omega^p(M))^G \rightarrow f\Gamma(M,\Omega^q(M))^G$ commutes with $\tilde{\Delta}$. Let
$\pi_p$ denote the projection of ${\bf H}^{0}_f(M,\Omega^p(M))^G$
onto $\mathrm{ker}\,\tilde{\Delta}$. By definition, on ${\bf
H}^{0}_f(M,\Omega^p(M))^G $,
\begin{equation}
\mathfrak{G}= \left(\tilde{\Delta}\Big|_{\mathrm{Im}\,
\tilde{\Delta}}\right)^{-1} \circ \pi_p.
\end{equation}
Now $T\tilde{\Delta}=\tilde{\Delta} T$ implies that
$T(\mathrm{ker}\,\tilde{\Delta})\subset
\mathrm{ker}\,\tilde{\Delta}$; and
$T(\mathrm{Im}\,\tilde{\Delta})\subset \mathrm{Im}\,\tilde{\Delta}$.
Hence

\begin{equation}
T\circ \pi_p = \pi_p \circ T.
\end{equation}
On the other hand,
\begin{eqnarray}
T\circ \left(\tilde{\Delta}\Big|_{\mathrm{Im}\,
\tilde{\Delta}}\right) &= & T\circ \tilde{\Delta}\circ (1-\pi_p)\cr
&=&  \tilde{\Delta}\circ  T\circ (1-\pi_p) = \tilde{\Delta}\circ
(1-\pi_p)\circ T \cr &=& \left(\tilde{\Delta}\Big|_{\mathrm{Im}\,
\tilde{\Delta}}\right)\circ T.
\end{eqnarray}
So on $\mathrm{Im}\, \tilde{\Delta}$,
\begin{equation}
T\circ \left(\tilde{\Delta}\Big|_{\mathrm{Im}\,
\tilde{\Delta}}\right)^{-1} =
\left(\tilde{\Delta}\Big|_{\mathrm{Im}\,
\tilde{\Delta}}\right)^{-1}\circ T.
\end{equation}
Therefore $\mathfrak {G}$ commutes with $T$.
\end{proof}

Finally we have

\begin{proposition}\label{prop:harmonic} Let $\mathfrak{H}^*(M)^G$ denote the kernel of the operator $\tilde{\Delta}$.
The map $H$ induces an isomorphism  $H:H^p(\Omega^*(M)^G,\, d)\to
\mathfrak{H}^*(M)^G$.
\end{proposition}

\begin{remark}
We remark that every element in $\mathfrak{H}^*(M)^G$ is of the form
$f\alpha$, where $\alpha$ is a $G$-invariant closed form.  For
$f\alpha\in \mathfrak{H}^*(M)^G$, as $\tilde{\Delta}(f\alpha)=0$,
\[
0=\langle\tilde{\Delta}(f\alpha), f\alpha\rangle_0=\langle d_f (f\alpha), d_f (f\alpha)\rangle_0+\langle d_f^*(f\alpha), d_f^*(f\alpha)\rangle_0.
\]
We conclude that $d_f(f\alpha)=fd\alpha=0$ and $d_f^*(f\alpha)=0$.
As $\alpha$ is $G$-invariant, we conclude that $d\alpha=0$.
\end{remark}

\begin{proof}
Let $\alpha$ be a $G$-invariant smooth closed $p$-form on $M$.
Consider $f\alpha\in {\bf H}_f^0(M, \Omega^p(M))^G$. As $d\alpha=0$,
$d_f (f\alpha)=0$. We have the following decomposition
(\ref{HodgeDecomp}),
\begin{equation}
f\alpha = d_f d_f^\ast \mathfrak {G}(f\alpha) +  d_f^\ast d_f
\mathfrak {G}(f\alpha) + H(f\alpha).
\end{equation}
Since $d_f$ commutes with $\tilde{\Delta}$, it commutes also with
$\mathfrak {G}$, so
\begin{equation}
f\alpha = d_f d_f^\ast \mathfrak {G}(f\alpha) +  d_f^\ast \mathfrak
{G} (d_f (f\alpha)) + H(f\alpha).
\end{equation}
Thus if $f\alpha$ is closed for $d_f$ (i.e., $d\alpha=0$), then
\begin{equation}\label{decomp}
f\alpha = d_f d_f^* \mathfrak {G}(f\alpha) + H(f\alpha).
\end{equation}
We define $H(\alpha)$ to be $H(f\alpha)$.

If $\alpha=d\beta$, then we have
\[
d_f (f\beta)=f\alpha.
\]
As the Green operator $\mathfrak {G}$ commutes with $d_f$, we have
$d_f d_f^* \mathfrak {G}(d_f(f\beta))=d_fd_f^*d_f(\mathfrak
{G}(f\beta))=d_f(d_f^*d_f+d_fd_f^*)\mathfrak
{G}(f\beta)=d_f(\tilde{\Delta}\mathfrak
{G}(f\beta))=d_f(f\beta-H(f\beta))$. Notice that elements in
$\ker\tilde{\Delta}$ are $d_f$-closed. So we have
$d_fd_f^*\mathfrak {G}(f\alpha)=d_f(f\beta)=f\alpha$, which shows
that $H(f\alpha)=0$.  This means that $H$ is a well-defined map from
$H^p(\Omega^*(M)^G,\, d)$ to $\mathfrak{H}^*(M)^G$.

If $H(f\alpha)=0$, then by Equation (\ref{decomp}), we have
$f\alpha=d_fd_f^*\mathfrak {G}(f\alpha)$. By Proposition
\ref{prop:reg}, we can write $d_f^*\mathfrak{G}(f\alpha)=f\beta$ for
a $G$-invariant smooth form $\beta$. Then $f\alpha=fd\beta$, and
$\alpha=d\beta$. This implies that $H$ is injective.

By the regularity property for $\tilde{\Delta}$ (Proposition
\ref{prop:reg}), elements in $\mathfrak{H}^*(M)^G$ are all smooth.
Furthermore, all elements in $\ker\tilde{\Delta}$ vanish under
$d_f$. So every element in $\mathfrak {H}^*(M)^G$ can be  written as
$f\alpha$, where $\alpha$ is a $G$-invariant smooth closed form. As
the image of $f\alpha\in \ker{\tilde{\Delta}}$ under $H$ is
$f\alpha$, we conclude that $H$ is onto.

\end{proof}

Theorem 1.1 is a corollary of Proposition \ref{prop:inv-hom},
\ref{prop:fredholm}, and \ref{prop:harmonic}.

\section{Euler Characteristic of a Proper Cocompact Action}

The finite dimensionality of the de Rham cohomology groups of
$G$-invariant differential forms allows us to define the \it Euler
characteristic \rm of such a proper cocompact $G$-action:

\begin{equation}\label{defchi}
\chi(M ;G) :=\sum_{i=0}^n (-1)^i \dim H^i(\Omega^*(M)^G,\, d).
\end{equation}

Our main result in this section is the following

\begin{theorem}
Let $M$ be a $n$-dimensional manifold on which a Lie group $G$ acts
properly and cocompactly.
\begin{enumerate}
\item[(i)]  We consider a family of twisted $G$-invariant de Rham differential operators
defined on $\Omega^*(M)^G$:
\begin{equation}\label{dtwist}
d_{A,k}(\alpha)=A^{-k} d(A^k\alpha)=(d+k A^{-1} dA\wedge)\alpha,\
\text{such that}\ d_{A,k}^2=0.
\end{equation}
Define the twisted cohomlogies of $G$-invariant differential forms
on $M$ to be the cohomology of the differential $d_{A,k}$, which is
denoted by $H^{p,k}(\Omega^*(M)^G; d)$. The cohomology group
$H^{p,k}(\Omega^*(M)^G;d)$ is finite dimensional.
\item[(ii)]
The Poincar\'e duality theorem holds for $k=1$, there is a
non-degenerate pairing between $H^{p,1}(\Omega^*(M)^G;d)$ and
$H^{n-p,1}(\Omega^*(M)^G;d)$. As a corollary, when $n$ is odd, the
Euler characteristic $\chi(M;G)$ of the proper cocompact $G$-action
on $M$ is 0;
\item[(iii)] When there is a nowhere vanishing $G$-invariant vector field on $M$,
the Euler characteristic of the $G$-action on $M$
is 0.
\end{enumerate}

\end{theorem}
We remark that by Proposition \ref{prop:projection},
$A^{-1}dA=d{\log(A)} $ is $G$-invariant. Therefore $d_{A,k}$ is
well-defined on $\Omega^*(M)^G$.

\begin{proof}

Our proof of statement (i) is a copy of the proof of Theorem
\ref{thm:main}. We define an operator $d_{f,A,k}$ on ${\bf H}_f^1$
generalizing (\ref{eq:dfn-df}) by
\[
d_{f,A,k}(f\alpha)=fA^{-k}d(A^{k}\alpha)=f(d+kA^{-1}dA\wedge)\alpha.
\]

We compute the adjoint of $d_{f,A,k}$. For $\alpha$ and $\beta$ two
$G$-invariant differential forms:

\begin{eqnarray*}
\langle d_{f,A,k}^* (f\alpha),\, f\beta \rangle_0 &=& \langle
f\alpha,\, d_{f,A,k} (f\beta) \rangle_0 = \langle f\alpha,\,
f A^{-k} d(A^k\beta) \rangle_0 \\
&=& \langle f^2 A^{-k} \alpha,\,  d(A^k\beta) \rangle_0 =\langle
\delta(f^2 A^{-k} \alpha),\,  A^k\beta \rangle_0\\
&=& \langle
f^2 A^{-k} \delta\alpha - A^{-k}2fi_{\nabla f}\alpha - f^2 (-k) A^{-k-1} i_{\nabla A}\alpha  ,\,  A^k\beta \rangle_0\\
&=& \langle
f  \delta\alpha - 2i_{\nabla f}\alpha +k f A^{-1} i_{\nabla A}\alpha  ,\,  f\beta \rangle_0\\
&=& \langle f  \delta\alpha - P_f(2i_{\nabla f}\alpha) + k f A^{-1}
i_{\nabla A}\alpha  ,\,  f\beta \rangle_0,
\end{eqnarray*}
and the $G$-invariance of $\alpha$ implies that
\begin{eqnarray}
P_f(2i_{\nabla f}\alpha)(x) & = & \frac{f(x)}{(A(x))^2}
\int_G\chi(g) f(gx) 2i_{\nabla f}g^*(\alpha)(x) dg\cr
 & =& \frac{f(x)}{(A(x))^2} \left(i_{\int_G\chi(g) 2 f(gx)\nabla f(gx) dg}
 \alpha\right)(x)\cr
  & =& \frac{f(x)}{(A(x))^2} \left(i_{\nabla\left(\int_G\chi(g)  (f(gx))^2  dg\right)}
 \alpha\right)(x)\cr
  & =& \frac{f(x)}{(A(x))^2} \left(i_{\nabla A^2}\alpha\right)(x)\cr
  &=& ( f A^{-2} i_{\nabla A^2}\alpha)(x)
  = 2 ( f A^{-1} i_{\nabla A}\alpha)(x).
\end{eqnarray}
Hence
\begin{equation}
d_{f,A,k}^* (f\alpha)= f(\delta+(k-2) A^{-1} i_{\nabla A} )\alpha.
\end{equation}

Now we define an operator from ${\bf H}_f^2$ to ${\bf H}_f^0$:
\begin{equation}
\widetilde{\Delta}_k=d_{f,A,k}d_{f,A,k}^*+d_{f,A,k}^*d_{f,A,k}.
\end{equation}
The analogues of Propositions
\ref{prop:fredholm}-\ref{prop:harmonic} for cohomology
$H^{p,k}(\Omega^*(M)^G;d)$ and $\widetilde{\Delta}_k$ in the Section
\ref{sec:hodge} easily generalize. Therefore any class in
$H^{p,k}(\Omega^*(M)^G; d)$ has a unique generalized harmonic form
representative $f\alpha$, i.e.
\begin{equation}
d_{f,A,k}(f\alpha)=d_{f,A,k}^*(f\alpha)=0.
\end{equation}
This proves that the dimension of $H^{p,k}(\Omega^*(M)^G;d)$ is
finite dimensional for any $p,k$.\\

For statement (ii), we prove that the pairing between ${\bf
H}_f^2(M, \Omega^p(M))^G$ and ${\bf H}_f^2(M, \Omega^{n-p}(M))^G$
induces a non-degenerate pairing on the space of generalized
harmonic forms of $\widetilde{\Delta}_1$, i.e.
\begin{equation}\label{eq:pairing}
(f\alpha, f\beta):=\int_M f\alpha\wedge f\beta.
\end{equation}
We prove that the Hodge star operator $*$ defines an isomorphism
between the space of generalized harmonic $p$-forms to the space of
generalized harmonic $(n-p)$-forms for the operator
$\widetilde{\Delta}_1$, which implies that the non-degeneracy of the
pairing (\ref{eq:pairing}).

We prove that $\widetilde{\Delta}_1\ast=\ast \widetilde{\Delta}_1$,
which implies that the Hodge star operator $\ast$ defines an
isomorphism between the generalized harmonic forms.

Using the following group of equations, where $\alpha$ is a
$G$-invariant $p$-form,
\begin{eqnarray*}
\delta \alpha &=& (-1)^{np+n+1} * d *\alpha,\\
* * \alpha & =& (-1)^{p(n-p)}\alpha,\\
i_{\nabla A} * ( *\alpha) &=&(-1)^{n-p} *(dA\wedge *\alpha),
\end{eqnarray*}
we can check that
\begin{eqnarray*}
d_f*(f\alpha)&=&(-1)^{p}*f\delta \alpha\\
f(\delta*\alpha)&=&(-1)^{p+1}*fd\alpha\\
fi_{\nabla A}*\alpha &=&(-1)^{p}*(dA\wedge \alpha)\\
fdA\wedge (*\alpha)&=& (-1)^{p+1}*(i_{\nabla A}\alpha).
\end{eqnarray*}
Combining the above equations, we have
\[
d_{f,A,1}(*f\alpha)=(-1)^{p}*d_{f,A,1}^*(f\alpha),\qquad\qquad
d^*_{f,A,1}(*f\alpha)=(-1)^{p+1}*d_{f,A, 1}(f\alpha).
\]
In particular, we have
\begin{equation}\label{eq:star-lap}
*\widetilde{\Delta}_1=\widetilde{\Delta}_1*.
\end{equation}

Equation (\ref{eq:star-lap}) shows that the Hodge star operator
commutes with the generalized Laplace operator
$\widetilde{\Delta}_1$. Therefore, the Hodge star operator $\ast$
defines an isomorphism between the kernels of $\widetilde{\Delta}_1$
on ${\bf H}_f^2(M, \Omega^{p}(M))^G$ and ${\bf H}_f^2(M,
\Omega^{n-p}(M))^G$.
Hence we have the Poincar\'e duality for the cohomology groups
$H^{p,1}(\Omega^*(M)^G; d)$:
\begin{equation}
H^{p,1}(\Omega^*(M)^G; d)\cong H^{n-p,1}(\Omega^*(M)^G; d)^*.
\end{equation}

For the statement about the Euler characteristic, we first notice
that $\chi(M ;G)$ is the index of the following Fredholm operator
\begin{equation}
d_f+d_f^*: {\bf H}_f^0(M, \Omega^{\mathrm{even}}(M))^G\rightarrow
{\bf H}_f^0(M, \Omega^{\mathrm{odd}}(M))^G.
\end{equation}

Since $\{d_{f,A,k}+d_{f,A,k}^*\}$ is a continuous family of Fredholm
operators with respect to $k\in{\bf R}$, their indices are all same.
This implies that
\begin{eqnarray}
\chi(M ;G)& =&\ind (d_{f,A,k}+d_{f,A,k}^*: {\bf H}_f^0(M,
\Omega^{\mathrm{even}}(M))^G\rightarrow {\bf H}_f^0(M,
\Omega^{\mathrm{odd}}(M))^G)\cr &=& \sum_{i=0}^n (-1)^i \dim
H^{p,k}(\Omega^*(M)^G,\, d).
\end{eqnarray}

When the dimension $n$ of the manifold $M$ is odd, the Poincar\'e
duality for $H^{p,1}(\Omega^*(M)^G; d)$ implies $\chi(M ;G)=0$.
Statement (ii) is thus proved.

\begin{remark}
When the Lie group $G$ is unimodular, then by replacing $f$ by
$\displaystyle \frac{f}{A}$, one can make all the differential
operators $d_{f,A,k}$ as well as the cohomology groups $H^{p,k}$
independent of $k$, therefore we have actually the Poincar\'e
duality for $H^{p}(\Omega^*(M)^G; d)$.
\end{remark}

Now statement (iii), we suppose that there is a nowhere vanishing
vector field $V$ on $M$. Without loss of generality, we may assume
that $|V|\equiv 1$ everywhere on $M$. We assume that $\{e_i\}$ is
an orthonormal basis of $TM$, and $\nabla^{TM}$ the Levi-Civita
connection of   the $G$-invariant Riemannian  metric.

Following \cite{A},  we calculate (cf. \cite[p. 73]{Z})
\begin{eqnarray}\label{eq:twist-D}
\hat{c}(V)(d_f+d_f^*)\hat{c}(V)(f\alpha) &=& f\left(
\hat{c}(V)(d+\delta)\hat{c}(V)\alpha - 2 A^{-1}\hat{c}(V) i_{\nabla
A} (\hat{c}(V)\alpha) \right)\cr &=& f\Big{(}-(d+\delta)\alpha+
\hat{c}(V)\sum_{i=1}^n c(e_i)\hat{c}(\nabla_{e_i}^{TM}V)\alpha
 \cr &  & \,\,\,\, -2 A^{-1} (V^*\wedge+ i_V)(V(A)\alpha -V^*\wedge
i_{\nabla A}\alpha+ i_{\nabla A}i_V\alpha)\Big{)}\cr &=&
f\Big{(}-(d+\delta)\alpha+ \hat{c}(V)\sum_{i=1}^n
c(e_i)\hat{c}(\nabla_{e_i}^{TM}V)\alpha
 \cr &  & \,\,\,\,+ 2 A^{-1}i_{\nabla A}\alpha -2 A^{-1}
 ( V(A)\hat{c}(V)\alpha+\hat{c}(V)i_{\nabla A}i_V\alpha)\Big{)}\cr
&=& -(d_f+d_f^*)(f\alpha)+ f\Big{(}\hat{c}(V)\sum_{i=1}^n
c(e_i)\hat{c}(\nabla_{e_i}^{TM}V)\alpha \cr &  & \hspace{4cm} -2
A^{-1}
 ( V(A)\hat{c}(V)\alpha+\hat{c}(V)i_{\nabla A}i_V\alpha)\Big{)}.
\end{eqnarray}
In the above formula, $c(v)$ and $\hat{c}(V)$ are the Clifford operators of
the vector field $V$ on the space $\Omega^{\text{odd}}$
and $\Omega^{\text{even}}$. More explicitly, if $V^*$ is the 1-form
dual to the vector field $V$ with respect to the riemannian metric,
$$c(V)(\alpha)=V^*\wedge \alpha-i_V\alpha,\qquad \hat{c}(V)(\alpha)=V^*\wedge \alpha+i_V\alpha,\qquad \alpha \in \Omega^*(M).$$

The above computation (\ref{eq:twist-D}) shows that the difference
between $\hat{c}(V)(d_f+d_f^*)\hat{c}(V)$ and $-(d_f+d_f^*)$ is an
operator of order 0. Proposition \ref{prop:fredholm} generalizes
directly to this operator and states that
$\hat{c}(V)(d_f+d_f^*)\hat{c}(V)$ is a Fredholm operator from ${\bf
H}_f^0(M, \Omega^{\text{odd}}(M))^G$ to ${\bf H}_f^0(M,
\Omega^{\text{even}}(M))^G$. Furthermore, as the difference between
$\hat{c}(V)(d_f+d_f^*)\hat{c}(V)$ and $-(d_f+d_f^*)$ has an order
less than the order of $-(d_f+d_f^*)$, one can also prove that the
operator
\[
-(d_f+d_f^*)+\epsilon
(\hat{c}(V)(d_f+d_f^*)\hat{c}(V)+(d_f+d_f^*)):{\bf H}_f^0(M,
\Omega^{\text{odd}}(M))^G\longrightarrow {\bf H}_f^0(M,
\Omega^{\text{even}}(M))^G
\]
is a Fredholm operator for any $\epsilon \in \reals$.

By the stability of index of Fredholm operators, we have
\begin{eqnarray}
\chi(M;G) & = & \ind (d_f+d_f^*: {\bf H}_f^0(M,
\Omega^{\mathrm{even}}(M))^G\rightarrow {\bf H}_f^0(M,
\Omega^{\mathrm{odd}}(M))^G)\cr & = & \ind
(\hat{c}(V)(d_f+d_f^*)\hat{c}(V): {\bf H}_f^0(M,
\Omega^{\mathrm{odd}}(M))^G\rightarrow {\bf H}_f^0(M,
\Omega^{\mathrm{even}}(M))^G)\cr & = & \ind
(-(d_f+d_f^*)+\text{lower order terms}:\cr &&\hspace{3cm} {\bf
H}_f^0(M, \Omega^{\mathrm{odd}}(M))^G\rightarrow {\bf H}_f^0(M,
\Omega^{\mathrm{even}}(M))^G)\cr & = & \ind (-(d_f+d_f^*): {\bf
H}_f^0(M, \Omega^{\mathrm{odd}}(M))^G\rightarrow {\bf H}_f^0(M,
\Omega^{\mathrm{even}}(M))^G)\cr & = & -\chi(M; G).
\end{eqnarray}
We conclude that $\chi(M;G)=0$.

\end{proof}

\vspace{2mm}

{\small \noindent{Xiang Tang}, Department of Mathematics, Washington
University, St. Louis, MO, 63130, U.S.A.,
Email:xtang@math.wustl.edu.

\vspace{2mm}

\noindent{Yi-jun Yao}, Mathematics Department, Pennsylvania State
University, State College, PA 16802, U.S.A.,
Email:yao@math.psu.edu., and School of Mathematical Sciences,
Fudan University, Shanghai 200433, P.R.China., Email:
yaoyijun@fudan.edu.cn.

\vspace{2mm}

\noindent{Weiping Zhang}, Chern Institute of Mathematics and LPMC,
Nankai University, Tianjin 300071, P.R. China.,
Email:weiping@nankai.edu.cn.

}
\end{document}